\documentclass[12pt]{article}

\usepackage{color}
\usepackage[centertags]{amsmath}
\usepackage{amsfonts}
\usepackage{amsthm}
\usepackage{newlfont}
\usepackage{amscd}
\usepackage{amsgen}
\usepackage{amssymb}
\usepackage{stmaryrd}
\usepackage{amssymb,amsmath}
\usepackage{amscd}
\usepackage{enumerate}

\newlength{\defbaselineskip} 
\theoremstyle{plain}
\newtheorem{thm}{Theorem}[section]
\newtheorem{cor}[thm]{Corollary}
\newtheorem{con}[thm]{Conjecture}
\newtheorem{df}[thm]{Definition}
\newtheorem{lema}[thm]{Lemma}
\newtheorem{obs}[thm]{Proposition}
\newtheorem{exm}[thm]{Example}
\newtheorem*{tthm}{Main Theorem}

\newtheorem{rem}[thm]{Remark}
\newtheorem{pr}{Program}
\theoremstyle{definition} 
\theoremstyle{definition} \newtheorem{ex}{Example}[section] %

 \numberwithin{equation}{section}
\def\p{\mathbb{P}}

\def\z{\mathbb{Z}}

\def\c{\mathbb{C}}

\DeclareMathOperator{\End}{End}

\def\p{\mathbb{P}}

\def\ob{\begin{obs}}
\def\kob{\end{obs}}
\def\dow{\begin{proof}}
\def\kdow{\end{proof}}

\def\tw{\begin{thm}}
\def\ktw{\end{thm}}
\def\hip{\begin{con}}
\def\khip{\end{con}}
\def\lem{\begin{lema}}
\def\klem{\end{lema}}
\def\ex{\begin{exm}}
\def\prog{\begin{pr}}
\def\kprog{\end{pr}}
\def\wn{\begin{cor}}
\def\kwn{\end{cor}}
\def\uwa{\begin{rem}}
\def\kuwa{\end{rem}}
\def\kex{\end{exm}}
\def\dfi{\begin{df}}
\def\kdfi{\end{df}}
\begin{document}
\title{Algebraic varieties representing group-based Markov processes on trees}
\author{Mateusz Micha\l ek\footnote{The author is supported by a grant of Polish MNiSzW (N N201 2653 33).}}
\maketitle
\begin{abstract}
  In this paper we complete the results of Sullivant and Sturmfels
  \cite{SS} proving that many of the algebraic group-based models for
  Markov processes on trees are pseudo-toric. We also show in which
  cases these varieties are normal. This is done by the generalization
  of the discrete Fourier transform approach. In the next step,
  following Sullivant and Sturmfels, we describe a fast algorithm
  finding a polytope associated to these algebraic models. However in
  our case we apply the notions of sockets and networks extending the
  work of Buczy\'nska and Wi\'sniewski \cite{BW} who introduced it for the
  binary case of the group $\z_2$.
\end{abstract}
\section*{Introduction}
In the recent years phytogenetic trees have been paid much attention
to because of the increased interest in evolutionary processes. They
describe the DNA sequence changes and therefore they are a powerful
tool that helps biologists to explain the evolution. They allow also
for a better understanding of the relations between the species based
on their DNA structure. The study of molecular phylogenetics is quite new but
already advanced.


Phylogenetic trees turned out to be very interesting also for
mathematicians who recognized their relationship with algebraic
varieties. Cavender and Felsenstein \cite{CF} as well as Lake \cite{L}
pioneered the work in algebraic phylogenetics in the late 80s by
introducing the invariants, that are polynomials describing a variety.

In 1993 Evans and Speed \cite{ES} observed that there is a natural
group action of $\z_2\times\z_2$ on the nucleobases $\{A,C,G,T\}$,
that defines models well-known from the biology. The method was
generalized in \cite{SSE} by proving that the obtained variety is
pseudo-toric\footnote{In most papers concerning phylogenetics authors
  refer to pseudo-toric varieties as toric ones. In this paper we are
  going to discuses normality, so we will distinguish between these
  two notions.}.  Their results were based
on the discrete Fourier transform and strongly
relied on the hypothesis that the group acting on the space of states
is abelian. Surprisingly they are valid also for some sub-models of
the primary models. 
For example, the resulting varieties
  were known to be pseudo-toric for 2-Kimura and Jukes-Cantor
models even before the publication \cite{SS}. Still, it has not been
specified exactly for which group-based models the variety is
pseudo-toric. Further information on algebraic methods in phylogenetics can be found
for example in \cite{SS} or \cite{PS} and the references
therein.

Many results in phylogenetic algebraic geometry, like the description of the generators in \cite{SS},
were based on the pseudo-toricness blanket
assumption on varieties. The arguments which imply that algebraic varieties related to submodels of general group-based models are pseudo-toric were commonly believed to work for any
group-based model, although the authors agreed \cite{privat} that the
assumption underlying their work may not be fulfilled and should have
been stated more precisely. For an example in which the arguments implying pseudo-toricness of varieties related to submodels of general group-based models do not apply see the
appendix. This problem seems to be crucial for the group-based
Markov processes on trees that do not necessarily arise from the
biology.

Although the inspirations for phylogenetic models come from biology,
they can be defined in a purely
mathematical setting. We believe that, although mathematicians are
often motivated by the empirical results, the main object of their
study are still general theories. The fact that they are applied very
often in other fields or in the real life is
rather the special feature of mathematics than the aim of
mathematicians. In case of phylogenetics the first models were
introduced by biologists. Then a mathematical setting was developed in
which biological models were covered. However other models (like a
binary model) turned out to be also very interesting. The study of
these objects resulted in many beautiful theorems and unexpected
connections - see \cite{Xu}, \cite{BW}, \cite{Man} and \cite{DK}. We
believe that the setting of $G$-models introduced in this paper is a
good foundation for developing the theory of group-based models. They
are sufficiently general to cover all the group-based models of
interest and still, as it has been shown, have all the expected
properties.
 
The main aim of this paper is to prove that a large class (including
well known 2-Kimura and Jukes Cantor models) of group-based algebraic
models are pseudo-toric. We will do this by using not only abelian
groups, but also arbitrary groups that have a normal, abelian subgroup.

The main theorem of this paper \ref{main} proves the pseudo-toric
result in a general setting:
\begin{tthm}
  Let $H$ be a normal, abelian subgroup of a group $G\subset
  S_n=Sym(A)$, that acts transitively and freely on a set $A$ of $n$
  elements. Let $\widehat W$ be the space of matrices invariant with
  respect to the action of $G$ and let $W$ be the vector space spanned
  freely by elements of $A$. Then the algebraic model
  for $(T,W,\widehat W)$ is pseudo-toric for any tree $T$.
\end{tthm}

Moreover in such a setting the main theorem 26 of \cite{SS} is
applicable.

Furthermore we investigate which varieties are normal and thereby
toric. For necessary information concerning toric varieties see
\cite{Ful}, \cite{oda} or \cite{Cox}. It turns out that many examples
arising from abelian groups (like 3-Kimura model) are toric, but for
arbitrary groups (the case of 2-Kimura) this is not true, even in the
case of the simplest $3$-leaf tree.
 
This paper does not demand any knowledge from nonmathematical
sciences. Readers interested solely in algebraic geometry can skip
remarks concerning biological or probabilistic setting, although the
intuition form Markov processes can be helpful.
 
The structure of the article looks as follows. In the first section we recall
basic facts about constructing varieties associated to the
phylogenetic trees. We will try not only to state precise definitions
but also to give some intuitive meaning of introduced objects. In the
second section we present results concerning group-based models in the abelian case. We also restate and conceptualize notions of sockets and networks from \cite{BW}. In section three we define $G$-models that are generalizations
of the abelian case. In a slightly different setting the idea of $G$-models can be found in \cite{WDB}. The family of $G$-models contains such
examples as 2-Kimura, Jukes-Cantor and 3-Kimura model. Our results
apply not only to biological phylogenetic models but also to many
other Markov processes. In section four we present efficient
algorithms used for constructing polytopes of the pseudo-toric
models. In section five we prove facts concerning normality of a model
and we show examples in which the obtained variety is not normal. To
prove normality we reduce any tree to the simplest 3-leaf tree checking it by the use of the Macaulay computer program. In the last
section we state some open questions.

\section*{Acknowledgements}
I would like to thank very much Jaros\l aw Wi\'sniewski for introducing me to
the subject and giving many useful ideas and helpful advices. I also
thank Bernd Sturmfels and Seth Sullivant for their inspiring article,
important information on group-based models and encouraging
words.

\section{Preliminaries}
In this section we will introduce a general setup for the rest of the paper. We remind basic facts about varieties associated to phylogenetic trees. From the point of view of Markov processes each point of our variety can be considered as a generalization of a probability distribution on the space of states of leaves of the tree. First we fix a rooted tree $T$. In many biological cases this tree encodes the history of the evolution of DNA. In our setting we assume that the root distribution is uniform. The root is thought of as a common ancestor, each edge corresponds to mutation and leaves correspond to species that mutated the last, and so are observable.
If we speak about directed edges of a tree $T$ we always assume that they are directed from the root. However it turns out that the variety associated to a given tree does not depend on the orientation of edges (so neither on the root) for any $G$-model, see remark \ref{niezalodorient}.

We fix a finite set $A$.
In the biological setting the set $A$ corresponds to the possible observed states on a given position in the DNA sequence. In biology the elements of $A$ are denoted by A, C, G, T. We will be considering more general Markov processes, so $A$ can be an arbitrary finite set.
\subsection{Algebraic setting of phylogenetic trees}

In biology we assume that we do not know the state of inner vertices of trees. That is why we associate to each node an $A$-valued random variable. Of course such variables give as a distribution - a set of nonnegative real numbers that sum up to 1. Following \cite{4aut} we generalize this setting. We work over the field of complex numbers $\c$ and we omit the condition of summing to one. In this setting in each vertex $v$ of the tree we obtain a complex linear space $W_v\cong W$ with the distinguished basis given by elements of $A$. Random variables described above correspond to the set of points with real, nonnegative entries that sum up to one - this is a simplex.

Apart from the states, for a specific model we also define the possible transition mechanism along edges (like a Jukes-Cantor model or Kimura models).
 In biology this part of the data depends on the model of the mutation mechanism that we choose. In general, for a Markov process it is defined by an association of a bi-stochastic matrix to each edge. The entries of this matrix correspond to relative probabilities of mutations along edges. In an algebraic setting (see also \cite{4aut}) the bi-stochastic condition is replaced by the condition that the sum of all coefficients in columns and rows are equal. For now, let us forget about this condition and we will see later that in our case it will be always satisfied. For the time being we associate to each edge from $v_1$ to $v_2$ a matrix, that is a linear map from $W_{v_1}$ to $W_{v_2}$ in the chosen basis. As in each vertex we have a distinguished basis, we can view the space of morphisms form $W_{v_1}$ to $W_{v_2}$ as $W_{v_1}\otimes W_{v_2}$, namely:
$$\End(W_{v_1},W_{v_2})\ni (f:W_{v_1}\rightarrow W_{v_2})\rightarrow\sum_{(g_1,g_2)\in G_{v_1}\times G_{v_2}}g_2^*(f(g_1))g_1\otimes g_2,$$
where $G_{v_i}$ is the distinguished basis of $W_{v_i}$ and $g_i^*$ is an element of the basis dual to $g_i$. In a given phylogenetic model we consider only matrices from some subspace $\widehat W$ of $W\otimes W$. Often we assume for example that the matrix is symmetric or as already mentioned bi-stochastic in an algebraic sense.

The rest of this section fixes some notation and can be found in \cite{BW}.
Now to each edge $e$ of a given rooted tree $T$ we associate a copy of $\widehat W$ denoted by $\widehat W_e$ (in fact, we can associate different vector spaces to different edges and all the arguments will work).
\dfi
We define the space of all possible states of a tree:
$$W_V=\bigotimes_{v\in V}W_v,$$ where $V$ is the set of vertices of a tree, the parameter space:
$$\widehat W_E=\bigotimes_{e\in E}\widehat W_e,$$
where $E$ is the set of edges and the space of the states of leaves:
$$W_L=\bigotimes_{l\in L} W_l,$$
where $L$ is the set of all leaves of the model.
\kdfi
Let $(v_1,v_2)$ be an edge that connects two adjacent vertices $v_1$ and $v_2$ (directed from $v_1$ to $v_2$).
We consider a linear map:
$${\widehat{\psi}} ':\bigotimes_{(v_1,v_2)\in E} (W_{v_1}\otimes W_{v_2})\rightarrow W_V.$$
Intuitively this map associates to a particular choice of matrices over all edges the "probability distribution" on the set of all possible states of the tree.
In order to define this map we have to define it on the basis. As $A$ can be considered as a basis of $W$ the basis of $$\bigotimes_{(v_1,v_2)\in E} (W_{v_1}\otimes W_{v_2})$$
can be chosen as follows:
$$\bigotimes_{e=(v_1,v_2)\in E} (a^{v_1}_e\otimes a^{v_2}_e),$$ where for each $e$ adjacent to $v$ the elements $a^v_e$ are elements of basis of $W_v$. Such an element corresponds to setting a matrix over an edge $e=(v_1,v_2)$ that gives all probabilities of mutation equal to $0$ apart from the probability of changing from $a^{v_1}_e$ to $a^{v_2}_e$ along edge $e$ equal to $1$. Now it is easy to set the corresponding probability on the set of states:
$$\widehat \psi'(\otimes_{e=(v_1,v_2)\in E}(a^{v_1}_e\otimes a^{v_2}_e))=\otimes_{v\in V}a^v,$$
if for any $v$ all $a^v_e$ are equal for all edges adjacent to $v$ (we define $a_v$ as the vector that they are all equal to) and $0$ otherwise.
\dfi
We define the map:
$$\widehat \psi:\widehat W_E\rightarrow W_V,$$
as the restriction of $\widehat\psi'$. This corresponds to the restriction of the possible space of matrices over any edge.
\kdfi
The closure of the induced rational map
$$\psi:\prod_{e\in E}\p(\widehat W_e)\dashrightarrow\p(W_V),$$
is called the complete projective geometric model.

Inspired by biological setting mathematicians consider not the complete projective geometric model, but a model that takes into account that we do not know the states of inner vertices of the tree $T$. We define a map:
$$\delta:W\rightarrow \c,$$
such that $\delta(a)=1$. In other words $\delta=\sum_{a\in A} a^*$. Next we define a map that hides the inner vertices.
\dfi
Let
$$\pi_L:W_V\rightarrow W_L,$$
be a map defined as $$\pi_L=(\otimes_{l\in L} id_{W_l})\otimes(\otimes_{v\in N}\delta_{W_v}),$$
where the map $\delta_{W_v}$ acts like $\delta$ on the vector space associated to a node $v$.
\kdfi
\uwa
\emph{The idea of the map $\pi_L$ is the following:
Consider points of $W_V$ as generalizations of probability distributions on the space of all possible states of the tree $T$ and points of $W_L$ as generalizations of probability distributions on the space of possible states of the leaves of $T$. Let $p\in W_V$. The image $\pi(p)$ corresponds to the distribution that to a given state $S$ of leaves associates the sum of all probabilities (induced by $p$) of states of the tree $T$ that coincide with $S$ on leaves. This remark can be read literally if we restrict the spaces $W_V$ and $W_L$ to simplices described at the beginning of the section, where the points have got a probabilistic meaning.}
\kuwa
If we compose the map $\widehat\psi$ with $\pi_L$ we obtain a map from $\widehat W_E$ to $W_L$.
\dfi
The map described above induces a rational map:
$$\prod_{e\in E}\p(\widehat W_e)\dashrightarrow\p(W_L).$$
The projective geometrical model, or just a model, denoted by $X(T)$ is the closure of the image of this map.
\kdfi

\section{Group models - abelian case}
Most of this section is well known. However we repeat some of the arguments in order to fix the notation and introduce some ideas that will be needed in following sections. We also generalize the notions of "sockets" and "networks" introduced in \cite{BW}, what enables us to generalize some of the results form $\z_2$ to arbitrary abelian groups.

Let $T$ be a rooted, trivalent tree. Let $A$ be a set of letters (possible states). In their famous paper \cite{ES} Evans and Speed recognized a natural action of an abelian group $G$ on $A$ in biological case, namely the group $G=\z_2\times \z_2$ acts on $\{A,C,G,T\}$ transitively and freely. Their method was described in a more general setting in \cite{SSE}. From now on we assume that we have a transitive and free action of an abelian group $G$ on $A$. Of course, if we define the vector space $W$ with the basis $A$, then the action of $G$ on $A$ extends to the action of $G$ on $W$.
\dfi
For $g\in G$ let $A_g$ be a matrix (equivalently a linear map) corresponding to the action of $g$ on $W$.
\kdfi

By choosing one element of the set $A$ and associating it to a neutral element of $G$ we may make an action preserving bijection between the elements of $A$ and $G$. An element associated to $a\in A$ will be denoted by $g_a$.
This allows us to find another basis of $W$, indexed by characters of $G$. This is done by the discrete Fourier transform.
\dfi
Let $\chi\in G^*$ be any character of the group $G$. We define a vector $w_\chi\in W$ by:
$$w_{\chi}=\sum_{a\in\ A}\chi(g_a)a.$$
\kdfi
One can prove that the elements $w_{\chi}$ form a basis of $W$.
Let us notice that although the choice of the bijection between $A$ and $H$ is not cannonical, the one dimensional spaces spanned by $w_\chi$ are. Changing the bijection just multiplies each vector $w_\chi$ by $\chi(g)$ for some $g\in G$.

The group structure distinguishes also naturally a specific model, namely the vector space $\widehat W$. The space $\widehat W\subset \End W$, but as we have a distinguished basis on $W$ made of the elements of $A$, we have an isomorphism $W^*\approx W$, so as already mentioned we may identify $\End W\approx W\otimes W$. We have a natural action of $G$ on $W\otimes W$ - the action of $g$ is just $g\otimes g$:
$$g(\sum\lambda a_1\otimes a_2)=\sum\lambda g(a_1)\otimes g(a_2).$$
Now, we may define $\widehat W$ as the set of fixed points of this action - see also \cite{WDB}.
\uwa\label{cotoW^}
\emph{In other words we take only such matrices that satisfy the following condition for any $g\in G$:}

If we permute the columns and rows of matrices with a permutation corresponding to $g$ then we obtain the same matrix.
\kuwa
\uwa\label{niezalodorient}
\emph{One can see that if $A\in\widehat W$, than $A^T\in\widehat W$. This means that if we consider a tree $T$ with two different orientations than the models are exactly the same. If a point is the image of some element of the parameter space with respect to a given orientation than it is also the image of an element of the parameter space with respect to the second orientation. We just have to transpose matrices that are associated to edges with different orientation.}
\kuwa
\ex
Let us consider the group $\z_2=(0,1)$. We obtain matrices of the form:
\[
\left[
\begin{array}{cccccccc}
a&b\\
b&a\\
\end{array}
\right]
\]
This is a binary model well studied in \cite{BW}.
\kex
\ex
In this example we consider the case of particular interest in biology. The natural action on $A=\{A,C,G,T\}$ is the action of the group $\z_2\times\z_2$. Let us consider the group $G=\z_2\times\z_2$. We obtain a model:
\[
\left[
\begin{array}{cccccccc}
a&b&c&d\\
b&a&d&c\\
c&d&a&b\\
d&c&b&a\\
\end{array}
\right],
\]
that is a 3-Kimura model.
\kex

One can also check that using the basis $(w_{\chi})_{\chi\in \widehat G}$ the elements of $\widehat W$ are exactly the diagonal maps. This means that we can define elements of $\widehat W$ also indexed by characters of the group $G$:
\dfi
$$l_{\chi}(w_{\chi'})=
\begin{cases}
w_{\chi} & \chi=\chi' \\
0 & \chi\neq\chi'
\end{cases}$$.
\kdfi
It follows that these elements are a basis of $\widehat W$. One can also check that using the basis $A$ we have:
\ob\label{lwstarych}
$$l_\chi (a_0)=\sum_{a\in A} \chi(g_{a_0}^{-1}g_a)a.$$
\kob
We see that the vectors $l_\chi$ are independent from the choice of the bijection between $A$ and $G$, as the element $g_{a_0}^{-1}g_a$ is a unique element of $G$ that sends $a_0$ to $a$, hence does not depend on the bijection. The map $l_\chi$ is a projection onto the (canonical) one dimensional subspace spanned by $w_\chi$.

Using this basis we will see that the map defining the model is given by a subsystem of the Segre system. First we need some lemmas. Of course, the action of $G$ on $W$ extends to the action of $G$ on $W_V$ and $W_L$. We have:
 \lem\label{wymiar}
$$\dim W_V^G=\frac{1}{|G|}|G|^{|V|}=|G|^{|V|-1}$$
$$\dim W_L^G=\frac{1}{|G|}|G|^{|L|}=|G|^{|L|-1}$$
\klem
\dow
Let us consider the following basis of $W_V$:
$$(\otimes_{v\in V} w_{\chi_v}).$$
The action of $g$ in this basis is diagonal, so the space of invariant vectors is spanned by invariant elements of this basis. As $g(w_\chi)=\chi(g^{-1}) w_\chi$ we obtain:
$$g(\otimes_{v\in V} w_{\chi_v})=\otimes_{v\in V} \chi_v(g^{-1})w_{\chi_v}=\prod_{v\in W}\chi_v(g^{-1})\otimes_{v\in V}w_{\chi_v},$$
so such an element is invariant iff for any $g\in G$ we have $\prod_{v\in W}\chi_v(g)=1$. This is equivalent to the condition that $\sum_{v\in V}\chi_v$ is equal to a trivial character (we use additive notation for the group $G^*$). From this we see that the dimension $\dim W_V^G$ is equal to the number of sequences indexed by vertices of the tree of characters that sum up to a neutral character. 
This gives us $|\widehat G|^{|V|-1}$ sequences what proves the first equality, as for abelian groups $|\widehat G|=|G|$.
The proof of the second equality is the same.
\kdow
\uwa
\emph{The basis $\{\otimes_{v\in V}w_{\chi_v}\}$ of $W_V$ in not natural (that is depends on the choice of the bijection between the set $A$ and $G$). However the basis $\{\otimes_{v\in V}w_{\chi_v}:\sum_{v\in V}\chi_v=\chi_0\}$ of $W_V^G$ is natural. Changing the bijection multiplies $w_\chi$ by $\chi(g)$ for a fixed $g\in G$, so as $\sum_{v\in V}\chi_v=\chi_0$, then of course $(\sum_{v\in V}\chi_v)(a)=1$ and the vectors remain unchanged.}
\kuwa

One can easily see that the image of $\widehat W_E$ in $W_V$ is invariant with respect to the action of $G$.
\ob\label{Segre}
The following morphism:
$$\psi:\prod_{e\in E}\mathbb{P}(\widehat W_e)\rightarrow \mathbb P(W_V^G),$$
is given by a full Segre system.
\kob
\dow
We prove that there is an isomorphism of $\widehat W_E$ and $W_V^G$. Thanks to lemma \ref{wymiar} the dimensions are the same, because the number of edges is equal to number of vertices minus one, so it is enough to prove the surjectivity of the dual morphism $W_V^*\rightarrow \widehat W_E^*$. In this proof we always consider the basis of $W$ given by elements of $A$. The space $\widehat W_E^*$ is generated by tensor products of functions that for a given matrix in $\widehat W$ return a given entry of the matrix:
$$\otimes_{e\in E} (a'_e)^*\otimes (a''_e)^*$$
where $a^*$ is an element dual to $a$ 
and if $e=(v_1,v_2)$, then $a'_e\in W_{v_1}$ and $a''_e\in W_{v_2}$ are elements of the chosen basis. Let $g_e$ be a unique element of $G$ that sends $a'_e$ to $a''_e$. 
Next, we construct an element of $W_V$ inductively starting from the root, where we take any element from $A$. Suppose that we have already defined an element $a$ for a vertex $v$. We take an edge $e=(v,w)$ and we define the state at $w$ as $g_ea$. One can easily check that the image of the dual of the tensor product of so defined elements gives us $\otimes_{e\in E} (a'_e)^*\otimes (a''_e)^*$, what proves the theorem.
\kdow
Moreover the isomorphism of $\widehat W_E$ and $W_V^G$ has got a very nice description in the new basis.
\ob\label{nowewspolpierwmorf}
The isomorphism of vector spaces $\widehat W_E$ and $W_V^G$ takes the base $\{\otimes_{e\in E} |G|l_{\chi_e}\}$ bijectively onto the base $\{\otimes_{v\in V}w_{\chi_v}:\sum_{v\in V}\chi_v=\chi_0\}$, where $\chi_0$ is a trivial character. This means that in these basis the Segre embedding is given by monomials.
\kob
\dow
We can see that:
$$l_{\chi}(a_0)=\frac{1}{|G|}\chi(g_{a_0}^{-1})w_{\chi}=\frac{1}{|G|}\sum_{a\in A}\chi(g_{a_0}^{-1}g_a)a.$$
This shows that the image of $\otimes_{e\in E} |G|l_{\chi_e}$ has got $c=\prod_{e=(v_1,v_2)\in E}\chi_e(g_{a_{v_1}}^{-1}g_{a_{v_2}})$ as a coefficient of $\otimes_{v\in V}a_v$. The coefficient $c$ is equal to: $$c=\prod_{e=(v_1,v_2)\in E}(-\chi_e)(g_{a_{v_1}})\chi_e(g_{a_{v_2}}).$$
For given characters $\chi_e$ let us define characters $\chi_v$ for all $v$ vertices of the tree as:
$$\chi_v=\sum_{(v,w)\in E}\chi_{(v,w)}-\sum_{(w,v)\in E}\chi_{(w,v)}.$$
This corresponds to summing all characters on edges adjacent to $v$ with appropriate signs, depending on the orientation of the edge. Using this notation:
$$c=\prod_{v\in V}\chi_v(g_{a_v}).$$

We define an element $\otimes_{v\in V}w_{\chi_v}$ that is clearly in the chosen basis of $W_V^G$ as each character $\chi_e$ is taken twice with different signs, so the sum of all $\chi_v$ is a trivial character. Moreover
$$\otimes_{v\in V}w_{\chi_v}=\otimes_{v\in V}(\sum_{a\in A}\chi_v(g_a)a),$$
so the coefficient of $\otimes_{v\in V}a_v$ is equal to $\prod_{v\in V}\chi_v(g_{a_v})$, what proves the theorem.
\kdow

Let us notice that apart from the action of $G$ on $W\otimes W$ given by $g\otimes g$ that allowed us to define $\widehat W$, we have got another action of $G$ on $W\otimes W$ given by $g\otimes id$, where $id$ is the identity map.
\lem
The action $g\otimes id$ restricts to $\widehat W$.
\klem
\dow
It is enough to prove that the image of the action of $g\otimes id$ on any element that is invariant with respect to the action $g'\otimes g'$ is also invariant. Let $C$ be any element of $\widehat W$.
$$(g'\otimes g')((g\otimes id)C=(g'g\otimes g')(C)=(gg'\otimes g')(C)=(g\otimes id)(g'\otimes g')(C)=$$$$(g\otimes id)(C).$$
Here we used the fact that $G$ is abelian.
\kdow
\dfi
Let $N$ be the set of nodes of the tree.
We define $\rho_v^g$ for $v$ an inner vertex and $g\in G$ as an isomorphism of the space $\widehat W_E$ that acts on edges $(v,w)$ as $g\otimes id$, on edges $(w,v)$ as $g^{-1}\otimes id$ and on edges not adjacent to $v$ as identity.
We also define a group $G_N$ as a group generated by $\rho_v^g$.
\kdfi
\uwa\label{jakgo1}
\emph{It is crucial to realize how $g\otimes id$ acts on elements of $\widehat W$ considered as morphisms. One can check that $g\otimes id (A_{g'})=A_{g'}\circ A_{g^{-1}},$
so the action of $g\otimes id$ composes the given morphism with $A_{g^{-1}}$.}
\kuwa
\uwa
\emph{The definition of $\rho_v^g$ corresponds to changing the matrices above edges adjacent to $v$ in such a way that they give the "probability" as if the state of $v$ was different. We will see in the following proposition that these definitions are useful, because if we do not know the states of inner vertices, then the orbits of the action of $G_N$ give us the same "probability distribution" on the space of states of leaves.}
\kuwa
First we need a technical lemma.
\lem
The group $G_N\cong G^{|N|}$. $G_N$ acts on $\widehat W_E$. There is a base in which $G_N$ acts diagonally.
\klem
\dow
We will prove that the good base is the tensor product of endomorphisms $l_\chi$. Using \ref{jakgo1} we obtain:
$$(g\otimes id (l_\chi))(w_{\chi'})=l_\chi A_{g^{-1}}(w_{\chi'})=$$
$$=l_\chi A_{g^{-1}}(\sum_{a\in A}\chi'(g_a)a)=l_\chi(\sum_{a\in A}\chi'(g_a)g^{-1}a)=$$
$$=l_\chi (\sum_{a\in A}\chi'(g_ag)a)=\chi'(g)l_\chi(w_{\chi'})=\chi(g) l_{\chi}(w_{\chi'}),$$
where the last equality follows from the fact that $l_\chi(w_{\chi'})$ is non zero only if $\chi=\chi'$. This proves that $g\otimes id (l_\chi)=\chi(g) l_{\chi}$, what proves the theorem.
\kdow
\uwa\label{bazakraw}
\emph{As described above the elements of the base of $\widehat W_E$ are bijective with the sequences of characters indexed by edges of a tree. In other words an element of a basis of $\widehat W_E$ can be described as an association of a character of $G$ to each edge of a tree. Moreover the elements of the basis of $\widehat W_E$ that are invariant with respect to the action of $G_N$ are exactly such associations that a (signed) product of characters around each inner vertex gives a trivial character. Precisely, for each vertex $v$ and edges $e_0=(w_0,v)$, $e_1=(v,w_1),\dots,e_i=(v,w_i),\dots,e_k=(v,w_k)$ we have:
$$-\chi_{e_0}+\sum_{i=1}^k\chi_{e_k}=\chi_0,$$
where the sum is the sum in $\widehat G$ (remember about the additive notation in $\widehat{G}$), $\chi_e$ is the character associated to an edge $e$ and $\chi_0$ is a trivial character.}
\kuwa
\ob\label{nowewspoldrugimorf}
The map:
$$W_V\ni\otimes_{v\in V}g_v\rightarrow \otimes_{l\in L}g_l\in W_L$$
in the basies $w_\chi$ is given by:
$$W_V\ni\otimes_{v\in V}w_{\chi_v}\rightarrow |G|^{|N|}\otimes_{l\in L}w_{\chi_l}$$
if all the characters for the inner vertices are trivial or zero otherwise.
\kob
\dow
First let us look at $\otimes_{v\in V}w_{\chi_v}$ in the old coordinates:
$$\otimes_{v\in V}w_{\chi_v}=\otimes_{v\in V}(\sum_{a\in A}{\chi_v}(g_a)a)=\sum_{(a_u)_{u\in V}\in A^V}(\prod_{v\in V}{\chi_v}(g_{a_v}))(\otimes_{v\in V} a_v),$$
where the sum $\sum_{(a_u)_{u\in V}\in A^V}$ is taken over all $|V|$-ples (indexed by vertices) of basis vectors. In other words this sum parameterizes the basis of $W_V$ made of tensor products of base vectors corresponding to elements of $G$.
This is of course equal to:
$$\sum_{(a_u)_{u\in N}\in A^N}\sum_{(a_l)_{l\in L}\in A^L}\prod_{v\in N}{\chi_v}(g_{a_v})\prod_{f\in L}{\chi_f}(g_{a_f})\otimes_{v\in N} a_v\otimes_{f\in L} a_f.$$
We see that the image in $W_L$ is equal to:
$$\sum_{(a_u)_{u\in N}\in A^N}\sum_{(a_l)_{l\in L}\in A^L}\prod_{v\in N}{\chi_v}(g_{a_v})\prod_{f\in L}{\chi_f}(g_{a_f})\otimes_{f\in L} a_f=$$
$$(\prod_{v\in N}(\sum_{g\in G}{\chi_v}(g)))\sum_{(g_l)_{l\in L}\in G^N}\prod_{f\in L}{\chi_f}(g_l)\otimes_{f\in L} a_f.$$
The product $(\prod_{u\in N}(\sum_{g\in G}{\chi_u}(g)))$ is equal to zero unless all characters $\chi_u$ for $u\in N$ are trivial. In the latter case the product is equal to $|G|^{|N|}$. Of course
$$\sum_{(g_l)_{l\in L}\in G^N}(\prod_{f\in L}{\chi_f}(g_l))(\otimes_{l\in L} g_l)=\otimes_{l\in L}w_{\chi_l},$$
what proves the proposition.

\kdow

The following theorem is a direct generalization to arbitrary abelian groups of theorem 2.12 from \cite{BW}.
\tw\label{iso}
The spaces $(W_L^G)$ and $(\widehat W_E)^{G_N}$ are isomorphic.
\ktw
\dow
One can prove it using dimension argument, but it is better to look how the basis are transformed.
The base of $(\widehat W_E)^{G_N}$ is given by $\otimes_{e\in E} |G|l_{\chi_e}$, where the signed sum of all characters at any vertex is trivial. This, thanks to \ref{nowewspolpierwmorf}, by the morphism $\widehat W_E\rightarrow W_V$ is transformed (bijectively) into an independent set $\otimes_{v\in V} w_{\chi_v}$, where characters for inner vertices are trivial and the sum of all characters is trivial. Using \ref{nowewspoldrugimorf} the image of this set via $W_V\rightarrow W_L$ gives us (bijectivly) $|G|^{|N|}\otimes_{l\in L} w_{\chi_l}$, where the characters sum up to a trivial character. The last set is the basis of $W_L^G$.
\kdow
This motivates the following definitions of sockets and networks.
\dfi
A socket will be an association of characters from $\widehat G$ to each leaf such that the sum of all these characters is a trivial character.

A network will be an association of characters from $\widehat G$ to each edge such that the (signed) sum of characters at each inner vertex gives the trivial character.
\kdfi

Let us generalize the results on sockets and networks from \cite{BW}.
One can see that there is a natural bijection of networks and sockets: a network is projected to leaves and a socket is extended to a network using the summing condition around inner vertices. Moreover each network determines naturally an element of the basis of $(\widehat W_E)^{G_N}$ and each socket an element of the basis of $W_L^G$. The isomorphism of \ref{iso} just uses this natural bijection. Of course as $G$ was abelian instead of associating characters we can associate group elements.

In case of the group $\z_2$ (\cite{BW}) sockets were even subsets of leaves. That was associating $1$ to chosen leaves and $0$ to the other leaves. The condition that the subset has got even number of elements is just the condition that the elements from the group sum up to a neutral element. We see that this definition is compatible. Networks were subsets of edges such that there was an even number chosen around each inner vertex - this is also a condition of summing up to a neutral element around each inner vertex.

Using the theorem \ref{iso} we know that the model $X$ is the closure of the image of the rational map:
$${\c^{|G||E|}}\rightarrow \p^{(|L|-1)\times |G|},$$
where the coordinates of the domain are indexed by pairs $(e,\chi)$ for $e$ an edge and $\chi\in G^*$. The coordinates of the codomain are indexed by sockets (or equivalently networks). Let $M$ be the lattice with the basis given by pairs $(e,\chi)_{e\in E,\chi\in G^*}$. Let $P$ be a polytope that is a subpolytope of a unit cube and whose points correspond to networks. It follows that the model is given by the polytope $P$.

To each pair $(e,\chi)$ where $e$ is an edge and $\chi$ is a character of $G$ we can associate a one parameter subgroup that is given as a morphism from $M$ to $\z$ and is a dual vector to the vector of the base of $M$ that is indexed by the pair $(e,\chi)$. In particular for each leaf $l$ and character $\chi\in G^*$ we obtain a one parameter subgroup $\lambda_l^\chi$. For each $t\in \c^*$ the action of $\lambda_l^\chi(t)$ on the model extends to the action on $\p^{(|L|-1)\times |G|}\supset X$. The weight of this action on the coordinate indexed by a socket $s$ is either $0$ or $1$ depending on whether the socket $s$ associates to the leaf $l$ character $\chi$ (in this case $1$) or not (in this case $0$).

\uwa
\emph{In \cite{BW} the authors considered only one one parameter subgroup for each leaf although their group had two elements. Notice however that in our notation for the group $\z_2$ the weights of the action of $\lambda_l^0$ are completely determined by the weights of the action of $\lambda_l^1$ - one weights are negations of the others. In our notation the authors considered only $\lambda_l^1$.}
\kuwa
\section{Nonabelian Case}
The setting of this section is sufficiently general to cover many Markov processes, in particular this will be a generalization of the results of the previous section. However the inspiration is a 2-Kimura model, that is a phylogenetic model in which the transition matrices are of the following type:
\[
\left[
\begin{array}{cccccccc}
a&b&c&b\\
b&a&b&c\\
c&b&a&b\\
b&c&b&a\\
\end{array}
\right].
\]

In this case, as in the previous section, we also have an abelian group $H=\z_2\times \z_2$ that acts on the basis $(A,C,G,T)$ of a four dimensional vector space $W$. As we have seen the fixed points of the action of $H$ on $W\otimes W$ define a 3-Kimura model. We may however define a larger group $G$, namely $D_8$ that contains $H$ as a normal subgroup and the action of $G$ on $W\otimes W$ defines the 2-Kimura model (compare with \cite{WDB}). That is why we consider the following setting.

Let $A$ be an $n$-element set of letters. Let $G$ be a subgroup of $S_n=Sym(A)$ (not necessarily abelian) acting on $A$. Suppose moreover that the group $G$ contains a normal, abelian subgroup $H$ and the action of $H$ on $A$ is transitive and free.
Elements of $A$ once again correspond to states of vertices of a phylogenetic tree $T$. We define $W$ to be a complex vector space spanned freely by elements of $A$, namely
$W=\oplus_{a\in A} \c_a,$
where $\c_a$ is a field of complex numbers corresponding to one dimensional vector space spanned by $a\in A$.

What changes is the definition of $\widehat W$. We define elements of $\widehat W$ as matrices fixed not only by the action of $H$, but by the whole action of $G$. We use the notation assuming that $End(W)\subset W\otimes W$.
\dfi
Let
$$\widehat W=\{\sum_{a_i,a_j\in A} \lambda_{a_i,a_j} a_i\otimes a_j: \lambda_{a_i,a_j}=\lambda_{g(a_i),g(a_j)}\forall g\in G\}.$$
\kdfi
\uwa
\emph{The characterization of $\widehat W$ from \ref{cotoW^} is still valid.}
\kuwa
\uwa
\emph{The situation of the previous section corresponds to $G=H$.}
\kuwa
\uwa\label{ostroznie}
\emph{As before by choosing some element of $A$ we may make a bijection between $A$ and $H$. However this time we have to be very careful. The action of $G$ on $A$ (as permutation) will not generally be the same as the action of $G$ on $H$ (as a group). We fix one such bijection. An element associated to $a\in A$ will be denoted by $h_a\in H$.}
\kuwa
We will often use the following easy observation:
\lem\label{onlyone}
Let $h\in H$ be such an element that as a permutation sends $a$ to $b$, where $a,b\in A$. Then $h=h_bh_a^{-1}$.
\klem
\dow
Of course both elements send $a$ to $b$, so because $H$ acts on $A$ transitively and freely, they have to be equal.
\kdow
\dfi
A $G$-model, will be a phylogenetic model associated to a phylogenetic tree $(T,W,\widehat W)$.
\kdfi
Our aim is to prove that also in this generalized setting we will obtain  pseudo-toric varieties.
We will proceed in four steps.
\begin{enumerate}
\item We introduce a general method for constructing endomorphisms of $W$ from complex functions on $H$. We prove that under certain conditions (namely a function should be constant on orbits of the conjugation action of $G$ on $H$) the obtained endomorphism is in $\widehat W$.
\item We prove that some sums (over the orbits of the action of $G$ on $\widehat H$) of characters of $H$ are functions that can define elements of $\widehat W$. We also notice that we obtain a set of independent vectors of $\widehat W$.
\item Using dimension arguments we prove that the set defined in step 2 is in fact a basis.
\item Finally, using theorems from section 3, we prove, using the new coordinates that our variety is pseudo-toric.
\end{enumerate}
\uwa
\emph{Of course we can define $\widehat W_H$ as a vector space of matrices fixed by the action of $H$. From the previous section we know that the closure of the image of the map:
$$\psi: \prod_{e\in E}\p(\widehat {(W_H)_e})\rightarrow W_L,$$
is a pseudo-toric variety. Moreover we also found the base in which the described morphism is given by monomials. Of course $\widehat W\subset \widehat W_H$, so our aim is to prove that the restriction of the previous map is also given by monomials in certain base. We will use the base on $\widehat W_H$ to define the base of $\widehat W$.}
\kuwa
\subsection{Correspondence between functions on $H$ and endomorphisms of $W$ - step 1}
We are going to define some endomorphisms of $W$.
\dfi
Let $f:H\rightarrow \c$ be {\it any} function. We define:
$$l_f=\sum_{a,b\in A} f(h_a^{-1}h_b) a\otimes b.$$
\kdfi
\uwa
\emph{Notice that because of \ref{lwstarych} the previous definition is consistent with the definition of $l_\chi$ for $\chi \in\widehat H$. Moreover the vector $l_f$ depends only on the function $f$ and not the bijection between $A$ and $H$, as $h_a^{-1}h_b$ is the only element from $H$ the sends $a$ to $b$.}
\kuwa

\ob\label{staladajedobry}
Let us consider an action of $G$ on $H$:
$$(g,h)\rightarrow ghg^{-1}.$$
If $f$ is constant on orbits of this action then $l_f\in\widehat W$.
\kob
\dow
Consider any element $g\in G$. We focus on two entries of the matrix $l_f$, namely $(a_1,b_1)$ and $(a_2,b_2)$, where
$$g(a_1)= a_2\text{    and   }g(b_1)= b_2.$$
These entries are from the definition of $l_f$ respectively $f(h_{a_1}^{-1}h_{b_1})$ and $f(h_{a_2}^{-1}h_{b_2})$ so because of remark \ref{cotoW^} we want to prove that:
$$f(h_{a_1}^{-1}h_{b_1})=f(h_{a_2}^{-1}h_{b_2}).$$
Consider an element $gh_{b_1}h_{a_1}^{-1}g^{-1}$. Clearly it is an element of $H$ (because $H$ was a normal subgroup of $G$) that sends $a_2$ to $b_2$. From lemma \ref{onlyone} we obtain:
$$gh_{b_1}h_{a_1}^{-1}g^{-1}=h_{b_2}h_{a_2}^{-1}.$$
Because $f$ was constant on orbits of the action of $G$ and $H$ was abelian we get:
$$f(h_{a_1}^{-1}h_{b_1})=f(h_{a_2}^{-1}h_{b_2}),$$
what completes the proof.
\kdow
\subsection{Appropriate functions on $H$ - step 2.}
In the abelian case we considered characters of $H$. As $G$ was equal to $H$, these functions were of course constant on (one element) orbits of the action of $G$ on $H$. In a general case it may happen that we do not have an equality
$$\chi(ghg^{-1})=\chi(h).$$
Of course this equality holds if a character of $H$ extends to a character of $G$, but this is not always the case. If we define the vectors $l_\chi$ for $\chi\in\widehat H$ they may not be in $\widehat W$. To obtain the vectors in $\widehat W$ we will sum up some characters to obtain functions that satisfy the condition of \ref{staladajedobry}. Consider the action of $G$ on $\widehat H$:
$$(g,\chi)(h)=\chi(ghg^{-1}).$$
Let $O$ be the set of orbits of this action. Of course elements of $O$ give a {\it partition} of $\widehat H$. Let us define for each element $o\in O$ a function $f_o:H\rightarrow\c$.
\dfi
Let
$$f_o=\sum_{\chi\in o}\chi.$$
Here we are summing characters as complex valued functions, not as characters, so this is a normal sum, not a product.
We obtain
$$l_{f_o}=\sum_{\chi\in o} l_\chi.$$
\kdfi

\ob
The function $f_o$ satisfies the conditions of \ref{staladajedobry} (is constant on orbits of the action of $G$ on $H$ by conjugation).
\kob
\dow
$$f_o(g'hg'^{-1})=(\sum_{\chi\in o} \chi(g'hg'^{-1})=$$
$$=\sum_{\chi\in o} (g',\chi)(h)=\sum_{\chi\in o} \chi(h)=f_o(h),$$
because the action of $g'$ is a permutation of the orbit $o$.
\kdow
\wn\label{independent}
The vectors $l_{f_o}$ for $o\in O$ are in $\widehat W$. Moreover, as $l_\chi$ formed a basis of $\widehat W_H$, and $l_{f_o}$ are sums over a \textsl{partition} of this basis, they are independent.
\kwn
\ob
Any function $f$ that is constant on orbits of $O$ is a linear combination of functions $f_o$.
\kob
\dow
Of course $f$ can be uniquely decomposed into sum of characters of $H$:
$$f=\sum_{\chi\in\widehat H} a_{\chi}\chi.$$
We have to prove that coefficients of $\chi$ in the same orbit are the same. Let $(g,\chi_1)=\chi_2$. We know that for any $h\in H$ we have
$$\sum_{\chi\in\widehat H} a_{\chi}\chi(h)=f(h)=f(ghg^{-1})=$$
$$\sum_{\chi\in\widehat H} a_{\chi}\chi(ghg^{-1})=\sum_{\chi\in\widehat H} a_{\chi}(g,\chi)(h).$$
From the orthogonality of characters we see that $a_{\chi_1}=a_{\chi_2}$ what completes the proof.
\kdow
\wn\label{zostalo}
The number of orbits in $O$ (and so the number of vectors $l_{f_o}$) is equal to the number of orbits of the action of $G$ on $H$:
$$(g,h)\rightarrow ghg^{-1}.$$
\kwn
\dow
This follows from comparing dimensions of spaces of complex functions on $H$ that are constant on orbits.
\kdow
\subsection{Dimension of $\widehat W$ - step 3}
We are going to prove that the dimension of $\widehat W$ is equal to the number of orbits $|O|$. First let us note that all coefficients of any matrix in $\widehat W$ (in the basis $A$) are determined by coefficients in the first row (this follows from the second section). We see that $\dim \widehat W$ is equal to the number of independent parameters in the first row. Let $e$ be a fixed element of the set $A$ (corresponding to the first row of the matrix). The action of $G$ imposes some conditions, namely the coefficient in the $e$-th row and $a$-th colomn and the coefficient in the $e$-th row and $b$-th colomn for $a,b\in A$ have to be equal iff there exists an element $g\in G$ such that:
$$g(e)= e\text{    and    }g(a)= b.$$
By $h_c$ we will denote the (unique) element of $H$ that sends $e$ to $c$.
\lem\label{step2}
The following conditions are equivalent:
\begin{enumerate}
\item $\exists_{g\in G}$ that sends $e$ to $e$ and $a$ to $b$,
\item the elements $h_a$ and $h_b$ are in the same orbit with respect to the action $(g,h)=ghg^{-1}$.
\end{enumerate}
\klem
\dow
Of course $h_a$ and $h_b$ are in the same orbit iff $h_a^{-1}$ and $h_b^{-1}$ are in the same orbit. For the proof we concentrate on the second variant.

1 $\Rightarrow$ 2: From \ref{onlyone} we know that $gh_a^{-1}g^{-1}=h_b^{-1}$, because both elements send $b$ to $e$.

1 $\Leftarrow$ 2: Suppose that $gh_a^{-1}g^{-1}=h_b^{-1}$. Let $g'=h_b^{-1}gh_{g^{-1}(b)}$. Of course $g'$ sends $e$ to $e$, but $g'=gh_{a}^{-1}h_{g^{-1}(b)}$, hence it also sends $a$ to $b$.
\kdow
\ob
The dimension of $\widehat W$ is equal to the number of orbits $|O|$.
\kob
\dow
Equal parameters in the first row of matrices in $\widehat W$ correspond bijectively to orbits of the action of $G$ on $H$ from \ref{step2} and remarks at the beginning of this subsection. This, along with \ref{zostalo}, finishes the proof.
\kdow
\wn
The elements $l_{f_o}$ for $o\in O$ form a basis of $\widehat W$.
\kwn
\subsection{$G$-models are pseudo-toric - step 4}
Let us define a basis on $\widehat W_e$ made of vectors $l_{f_o}$. If we consider the inclusion map:
$$\widehat W_e\rightarrow \widehat {(W_H)}_e,$$
in the basis made respectively of $l_{f_o}$ and $l_{\chi}$, than this is just copying some coordinates, because each $l_{f_o}=\sum l_{\chi}$ and $o$ determines {\it uniquely} all $\chi$. Namely the coordinate of
$l_{f_{o_e}},$
is copied to all coordinates
$l_{\chi_e},$
for $\chi_e\in o_e$. One can also say that the coordinate of
$l_{\chi_e}$
is equal to the coordinate of
$l_{f_{o_e}},$
where $o_e$ is the orbit of $\chi_e$.
Of course this shows that the map from $\prod_{e\in E}\p(\widehat W_e)$ to $\p(W_L)$ that parameterizes the model is also given by monomials - these are exactly monomials from section 2, where we just make some variables equal to each other. We may also look at the following commutative diagram:
\[\begin{array}{lllll}
\prod_{e\in E}\p(\widehat W_e)&\rightarrow&\p(\widehat W_E)&\dashrightarrow&\p(W_L)\\
\hskip 30pt\downarrow&&\downarrow&&\updownarrow\\
\prod_{e\in E}\p(\widehat {W_H}_e)&\rightarrow&\p(\widehat {W_H}_E)&\dashrightarrow&\p(W_L)\\
\end{array}\]
This proves the main theorem of this section:
\tw\label{main}
Let $H$ be a normal, abelian subgroup of a group $G\subset S_n$, that acts transitively and freely on a set $A$ of $n$ elements. Let $\widehat W$ be the space of matrices invariant with respect to the action of $G$ and let $W$ be the vector space spanned freely by elements of $A$. Then a $G$-model of $(T,W,\widehat W)$ is pseudo-toric for any tree $T$.
\ktw

\section{Polytopes of $G$-models}
In this section we will show how the construction from the previous section works on Kimura models and we will present the algorithm for constructing a polytope of a model for a given group $G$ with a normal subgroup $H$. The method was described in a different language in \cite{SS}. The main difference (apart from the notation) is that the authors assumed the existence of a friendly labelling function, that described which characters are identified. In case of $G$-models we know precisely what this function is: it associates to a given character its orbit of a $G$ action. It is easy to see that this is a friendly labelling.

If $G=H$ this is particulary easy. The polytope has got $|G|^{|E|-|N|}$ vertices and the algorithm works in time $O(|N|(|G|^{|E|-|N|}))$ assuming that we can perform group operations in unit time.
\prog
\begin{enumerate}
\item Orient the edges of the tree from the root.
\item For each inner vertex choose one outcomming edge.
\item Make a bijection $b:G\rightarrow B\subset \z^{|G|}$, where $B$ is the standard basis of $\z^{|G|}$.
\item Consider all possible associations of elements of $G$ with not-chosen edges (there are $|G|^{|E|-|N|}$ such associations).
\item For each such associations, make a full association by assigning an element of $G$ to each chosen edge in such a way that the (signed) sum of elements around each inner vertex gives a neutral element in $G$.
\item For each full association output the vertex of the polytope: $(b(g_e)_{e\in E})$, where $g_e$ is the element of the group associated to edge $e$.
\end{enumerate}
\kprog
\ex\label{3Kim}
For a 3-Kimura model on a tree with one inner vertex
the vertices of $P$ correspond to triples of characters of the group that sum up to a neutral character:

1) $(0,0),(0,0),(0,0)$ \hskip30pt 2) $(0,0),(1,0),(1,0)$ \hskip30pt 3) $(1,0),(0,0),(1,0)$

4) $(1,0),(1,0),(0,0)$\hskip30pt  5) $(0,0),(0,1),(0,1)$\hskip30pt  6) $(0,1),(0,0),(0,1)$

7) $(0,1), (0,1),(0,0)$\hskip30pt  8) $(0,0),(1,1),(1,1)$\hskip30pt  9) $(1,1),(0,0),(1,1)$

10) $(1,1),(1,1),(0,0)$\hskip24pt  11) $(0,1),(1,0),(1,1)$\hskip24pt  12) $(0,1),(1,1),(1,0)$

13) $(1,0),(1,1),(0,1)$\hskip24pt  14) $(1,0),(0,1),(1,1)$\hskip24pt  15) $(1,1),(0,1),(1,0)$

16) $(1,1),(1,0),(0,1)$

This in the coordinates of the lattice gives us vertices of the polytope:

1) 1,0,0,0,1,0,0,0,1,0,0,0\hskip30pt 2) 1,0,0,0,0,1,0,0,0,1,0,0

3) 0,1,0,0,1,0,0,0,0,1,0,0\hskip30pt 4) 0,1,0,0,0,1,0,0,1,0,0,0

5) 1,0,0,0,0,0,1,0,0,0,1,0\hskip30pt 6) 0,0,1,0,1,0,0,0,0,0,1,0

7) 0,0,1,0,0,0,1,0,1,0,0,0\hskip30pt 8) 1,0,0,0,0,0,0,1,0,0,0,1

9) 0,0,0,1,1,0,0,0,0,0,0,1\hskip30pt 10) 0,0,0,1,0,0,0,1,1,0,0,0

11) 0,0,1,0,0,1,0,0,0,0,0,1\hskip25pt 12) 0,0,1,0,0,0,0,1,0,1,0,0

13) 0,1,0,0,0,0,0,1,0,0,1,0\hskip25pt 14) 0,1,0,0,0,0,1,0,0,0,0,1

15) 0,0,0,1,0,0,1,0,0,1,0,0\hskip25pt 16) 0,0,0,1,0,1,0,0,0,0,1,0

\kex
The basis for $\widehat W$ for 3-Kimura (in previous notation vectors $l_\chi=\sum\chi(h_a^{-1}h_b)a\otimes b$) is the following:
\[
\left[
\begin{array}{cccccccc}
1&1&1&1\\
1&1&1&1\\
1&1&1&1\\
1&1&1&1\\
\end{array}
\right],
\left[
\begin{array}{cccccccc}
1&-1&1&-1\\
-1&1&-1&1\\
1&-1&1&-1\\
-1&1&-1&1\\
\end{array}
\right],\]\[
\left[
\begin{array}{cccccccc}
1&-1&-1&1\\
-1&1&1&-1\\
-1&1&1&-1\\
1&-1&-1&1\\
\end{array}
\right],
\left[
\begin{array}{cccccccc}
1&1&-1&-1\\
1&1&-1&-1\\
-1&-1&1&1\\
-1&-1&1&1\\
\end{array}
\right].
\]
For the 2-Kimura model the elements of $H$ are in the order:
$$(1)(2)(3)(4);(1,2)(3,4);(1,3)(2,4);(1,4)(2,3)$$
and $G$ is spanned by $H$ and $(3,4)$.

Now if we consider the action of $G$ on $\widehat H$ we obtain:
\begin{itemize}
\item The orbit of the trivial character is of course only the trivial character. This tells us that the first vector is in $\widehat W_G$ and will be considered as the first basis vector.
\item The orbit of the character that associates $-1$ to $(1,3)(2,4)$ and $(1,4)(2,3)$ has got also only one element. For example let us notice that  $$\chi((3,4)(1,3)(2,4)(3,4))=\chi((1,4)(2,3))=-1=\chi((1,3)(2,4)).$$ This means that the last vector also will be a basis vector of $\widehat W_G$.
\item There is also one more orbit that contains two left characters. If we take their sum (as functions, not characters) we obtain a function that associates $2$ to $(1)(2)(3)(4)$, $-2$ to $(1,2)(3,4)$ and $0$ to other two elements. This gives us an element:
\[
\left[
\begin{array}{cccccccc}
2&-2&0&0\\
-2&2&0&0\\
0&0&2&-2\\
0&0&-2&2\\
\end{array}
\right]
\]
This is of course the sum of two $l_\chi$.
\end{itemize}
We obtain $f_1=l_1$, $f_2=l_2+l_3$, $f_3=l_4$, where $l_i$ are matrices introduced above. Let $F=\{f_1,f_2,f_3\}$ and $L=\{l_1,\dots,l_4\}$. From the previous section $F$ is the basis of $\widehat W_G$ and $L$ of $\widehat W_H$. This can be checked directly in this example.
Let us now look at the map for a tree $Y$ with one inner vertex. Elements of $\widehat W_G$ are special elements of $\widehat W_H$. We have a map:
$$(f^{e_i}_j)_{j=1,\dots,3,i=1,\dots,3}\rightarrow (l^{e_i}_j)_{j=1,\dots,4,i=1,\dots, 3}.$$
Here $j$ parameterizes base vectors and $i$ parameterizes edges.
Our model is the composition of this map and a model map for $H$.
The image of the first map is a subspace given by a condition that the coordinates corresponding to $l^{e_i}_2$ and $l^{e_i}_3$ are equal for each $i=1,\dots,3$.
Let us see this directly:

The fixed bijection $b$ is the following:
$$b(e)=(1,0,0,0),\quad b(\chi_3)=(0,1,0,0)$$
$$b(\chi_1)=(0,0,1,0),\quad b(\chi_2)=(0,0,0,1)$$
where $\chi_1$ and $\chi_3$ are in the same orbit. Now the domain of $\psi$ for the group $H$ is
$\{(x_1,\dots,x_{12}):x_i\in \c\}$ in the order corresponding to \ref{3Kim} (we fix an isomorphism with $\chi_1=(1,0)$ and $\chi_3=(0,1)$). This tells us that  the subspace $\prod_{e\in E} (\widehat W_G)_e$ is given by conditions $x_2=x_3$ (the coordinates of $l_2$ and $l_3$ for $\widehat W_H^{e_1}$), $x_6=x_7$, $x_{10}=x_{11}$.

This procedure works generally. After having fixed the polytope for a subgroup $H$, that is in the lattice $M$ (whose coordinates are indexed by edges and characters of $H$) we consider a morphism form $M$ onto the lattice $M'$ (whose coordinates are indexed by edges and orbits of characters of $H$) that just assigns a character to a given orbit. This morphism if of course just summing up coordinates that are in the same orbit of the action of $G$ on $\widehat H$. The image of the polytope $P$ is a polytope of our model. For 3-Kimura we sum up coordinates (order from \ref{3Kim}) obtaining a polytope for 2-Kimura model:

1) 1,0,0,1,0,0,1,0,0\hskip30pt 2) 1,0,0,0,1,0,0,1,0

3) 0,1,0,1,0,0,0,1,0\hskip30pt 4) 0,1,0,0,1,0,1,0,0

5) 1,0,0,0,1,0,0,1,0\hskip30pt 6) 0,1,0,1,0,0,0,1,0

7) 0,1,0,0,1,0,1,0,0\hskip30pt 8) 1,0,0,0,0,1,0,0,1

9) 0,0,1,1,0,0,0,0,1\hskip30pt 10) 0,0,1,0,0,1,1,0,0

11) 0,1,0,0,1,0,0,0,1\hskip25pt 12) 0,1,0,0,0,1,0,1,0

13) 0,1,0,0,0,1,0,1,0\hskip25pt 14) 0,1,0,0,1,0,0,0,1

15) 0,0,1,0,1,0,0,1,0\hskip25pt 16) 0,0,1,0,1,0,0,1,0

what after removing double entries gives vertices:

1) 1,0,0,1,0,0,1,0,0\hskip30pt 2) 1,0,0,0,1,0,0,1,0

3) 0,1,0,1,0,0,0,1,0\hskip30pt 4) 0,1,0,0,1,0,1,0,0

5) 1,0,0,0,0,1,0,0,1\hskip30pt 6) 0,0,1,1,0,0,0,0,1

7) 0,0,1,0,0,1,1,0,0\hskip30pt 8) 0,1,0,0,1,0,0,0,1

9) 0,1,0,0,0,1,0,1,0\hskip30pt 10) 0,0,1,0,1,0,0,1,0

\section{Normality of $G$-models}
We have seen that the models associated to a group containing a normal, abelian case are pseudo-toric. To obtain results that some varieties are toric one needs to prove normality. We will see that in general one cannot expect a $G$-model to be normal, but in many cases it is. The normality of varieties is equivalent to normality of polytopes. First let us start with a technical lemma.
\lem\label{prod}
Let $P_1$ and $P_2$ be two normal polytopes contained respectively in lattices $L_1$ and $L_2$ spanned by the points of the polytopes. Suppose that we have got morphisms $p_i:L_i\rightarrow L$ of lattices for $i=1,2$ such that $p_i(P_i)\subset S$, where $S$ is a standard symplex (convex hull of standard basis). Then the fiber product $P_1\times_{L}P_2$ is normal in the lattice spanned by its points.
\klem
\dow
Let $q\in n (P_1\times_{L}P_2)$ and let $q_i$ be the projection of $q$ to $L_i$. One can see that as $q$ was the sum of points that belong to $P_1\times_{L}P_2$ with coefficients summing up to $n$ and was in the convex hull of $n$ times the points of $P_1\times_{L}P_2$, then each $q_i$ is the sum of points that belong to $P_i$ with coefficients summing up to $n$ and is in the convex hull of $n$ times the points of $P_i$. This means that $q_i\in nP_i$. From the assumptions we obtain:
$$q_i=\sum_{j=1}^n v_j^i,$$
with each $v_j^i\in P_i$. We also know that $p_1(q_1)=p_2(q_2)$ and this is an element of $nS$. Moreover $p_i(v_j^i)\in S$. But let us notice that each element of $nS$ can be {\it uniquely} written as the sum of $n$ elements of $S$. This means that the collections $(p_1(v_1^1),\dots,p_1(v_n^1))$ and $(p_2(v_1^2),\dots,p_2(v_n^2))$ are the same up to permutation, so we can assume that $p_1(v_j^1)=p_2(v_j^2)$. So we can lift each pair $(v_j^1,v_j^2)$ to a point $v_j\in P_1\times_{L}P_2$ that projects respectively to $v_j^1$ and $v_j^2$. One obtains $q=\sum_{j=1}^n v_j$ what finishes the proof.
\kdow

Let us consider two trees $T_1$ and $T_2$ (in an abelian case) and two leaves $l_i\in T_i$, $i=1,2$. Let $P_i\in M_i$ be a polytope for a tree $T_i$ for $i=1,2$. Now we can take the lattice $L$ whose basis is given by characters of $G$. We take the projection $p_i:M_i\rightarrow L$ onto the coordinates indexed by pairs $(l_i,\chi)_{\chi\in G^*}$. From the construction of the polytope of the model we know that its vertices correspond to networks. The fiber product corresponds to such pairs of networks that associate the same character to $l_1$ and $l_2$. This is the same as associating a character to each edge of the tree $T$, that is obtained from $T_1$ and $T_2$ by gluing chosen leaves. We see that the fiber product of polytopes gives a polytope of a new tree.

This is also true in a general (not necessarily abelian) case. Now the vertices correspond to associations of orbits of characters of $H$ to each edge, such that we can find for each edge a representative and the representatives form a network. Let us notice that by the conjugation action of $G$ on $H^*$ we can always choose a representative on one fixed leaf arbitrarily. This means that if the orbits associated to leaves $l_1$ and $l_2$ agree than, we can find a network on the tree $T$ that projects to a fixed element of the orbit. This proves that also in a general case the polytope of the tree $T$ is given as a fiber product of polytopes for trees $T_i$.

\uwa
\emph{Of course one can notice  that the above assumptions are equivalent to the fact that the labelling function from \cite{SS} was friendly.}
\kuwa

The projections $p_i$ of the polytopes $P_i$ are indeed in $S$. This means that thanks to lemma \ref{prod} we only have to consider normality of a model for a 3-leaf tree, as each trivalent tree can be constructed by gluing such trees.
\ob
The $G$-models for abelian groups: $\z_2$, $\z_2\times \z_2$, $\z_3$ and $\z_4$ are normal (for trivalent trees).
\kob
\dow
One can easily find the polytopes for these models with the 3-leaf tree using the algorithm from the previous section. One can find the coordinates in the lattice that is bigger than the lattice spanned by the points of the polytope, but after an easy change of coordinates one can find the right presentation of the polytope. We can also compute the Gr\"obner basis of its cone and see that it is in the polytope (we use computer program Macaulay) what proves normality.
\kdow
\uwa
The polytope of a 2-Kimura model is not normal.
\kuwa
\dow
 The polytope of a 2-Kimura model is:

1) 1,0,0,1,0,0,1,0,0\hskip30pt 2) 1,0,0,0,1,0,0,1,0

3) 0,1,0,1,0,0,0,1,0\hskip30pt 4) 0,1,0,0,1,0,1,0,0

5) 1,0,0,0,0,1,0,0,1\hskip30pt 6) 0,0,1,1,0,0,0,0,1

7) 0,0,1,0,0,1,1,0,0\hskip30pt 8) 0,1,0,0,1,0,0,0,1

9) 0,1,0,0,0,1,0,1,0\hskip30pt 10) 0,0,1,0,1,0,0,1,0

 The point
$(1,0,1,1,0,1,1,0,1)$ is in the lattice spanned by the vertices of the polytope and in $2P$ but it is not sum of two points of $P$.
\kdow
\section{Open Problems}
From the previous section we see that all checked $G$-models for an abelian group $G$ are toric. It would be interesting to characterize groups for which $G$-models are toric. It is natural to ask whether each $G$-model for an abelian group is normal.
\section*{Appendix}
Here we show an explicit example when the equality of the parameters before the Fourier transform does not imply the equality after it.

Let us consider the group $G=\z^4$.
The space of matrices is of the form:
\[
\left[
\begin{array}{cccccccc}
a&b&c&d\\
d&a&b&c\\
c&d&a&b\\
b&c&d&a\\
\end{array}
\right].
\]
The matrix of the type above corresponds to a function $f:G\rightarrow\c$, such that $f(0)=a$, $f(1)=b$, $f(2)=c$ and $f(3)=d$.
Now the Fourier transform of $f$ gives us:
$\widehat f(\chi_0)=a+b+c+d$, $\widehat f(\chi_1)=a+ib-c-id$, $\widehat f(\chi_2)=a-b+c-d$, $\widehat f(\chi_3)=a-ib-c+id$. If we consider a submodel defined by $f(1)=f(2)$ what corresponds to $b=c$ the Fourier transform gives us respectively $(x_0,x_1,x_2,x_3)=(a+2b+d,a+(i-1)b-id,a-d,a-(i+1)b+id)$, what is a linear subspace defined as $(i+1)x_1-2ix_2+(i-1)x_3=0$, hence is not an equality of two distinct variables.

Mateusz Micha\l ek

{Mathematical Institute of the Polish Academy of Sciences,

\'{S}w. Tomasza 30, 31-027 Krak\'{o}w, Poland}
\vskip 2pt
{Institut Fourier, Universite Joseph Fourier,

100 rue des Maths, BP 74, 38402 St Martin d'H\`eres, France}

e-mail address:\emph{wajcha2@poczta.onet.pl}

\begin{thebibliography}{Wi2'}
\bibitem[BDW]{WDB} W. Buczy\'nska, M. Donten, J. Wi\'sniewski, Isotropic models of evolution with symmetries, \emph{Interactions of Classical and Numerical Algebraic Geometry: A Conference in Honor of A. J. Sommese}, May 22-24, (2008).
\bibitem[BW]{BW} W. Buczy\'nska, J. Wi\'sniewski, On the geometry of binary symmetric models of phylogenetic trees, \emph{J. European Math. Soc.} 9,609-635, (2007).
\bibitem[CF]{CF} J. Cavender, J. Felsenstein,  Invariants of phylogenies: A simple case with discrete states, \emph{J. Classif.} 4, 57-71, (1987).
\bibitem[CLS]{Cox} D. Cox, J. Little, H. Schenck, \emph{Toric varieties}, preprint at http://www.cs.amherst.edu/~dac/toric.html.
\bibitem[DK]{DK} J. Draisma, J. Kuttler, On the ideals of equivariant tree models. (English summary) Math. Ann. 344, no. 3, 619--644, (2009).
\bibitem[ERSS]{4aut} N. Eriksson, K. Ranestad, B. Sturmfels, S. Sullivant, Phylogenetic Algebraic Geometry,  \emph{Projective Varieties with Unexpected Properties}; Siena, Italy, 237-256, Berlin,   de Gruyter, (2004).
\bibitem[ES]{ES} S. Evans, T. Speed, Invariants of some probability models used in phylogenetic  inference, \emph{Annals of Statistics} 21, 355-377, (1993).
\bibitem[Ful]{Ful} W. Fulton, \emph{Introduction to toric varieties}, Princeton University Press, (1993).
\bibitem[L]{L} J. Lake, A rate-independent technique for analysis of nucleaic acid sequences: Evolutionary parsimony, \emph{Mol. Biol. Evol.} 4, 167-191, (1987).
\bibitem[Man]{Man} C. Manon, The algebra of Conformal Blocks, preprint at arXiv:0910.0577
\bibitem[Oda]{oda} T. Oda, \emph{Convex bodies and algebraic geometry}, Springer, (1987).
\bibitem[PS]{PS} L. Pachter, B. Sturmfels, \emph{Algebraic statistics for computational biology}, Cambridge University Press, (2005).
\bibitem[Ph]{Ph} C. Semple, M. Steel, \emph{Phylogenetics}, Oxford University Press, (2003).
\bibitem[Pr]{privat} B. Sturmfels, S. Sullivant, privat communication (2010).
\bibitem[SS]{SS} B. Sturmfels, S. Sullivant, Toric ideals of phylogenetic invariants, \emph{J. of Comp. Biol.} 12(2), 204-228, (2005).
\bibitem[SSE]{SSE} L. Szekely, P. Erd\"os, M. Steel, Fourier Calculus on Evolutionary Trees, \emph{Adv. in App.  Math.} 14, 200-216, (1993).
\bibitem[SX]{Xu} B. Sturmfels, Z. Xu, Sagbi Bases of Cox-Nagata Rings, preprint at http://arxiv.org/abs/0803.0892.
\end{thebibliography}
\end{document}